\documentclass[12pt]{amsart}
\usepackage{amsmath,amsfonts,amssymb}

\theoremstyle{plain}

\newcommand{\suml}{\mathop{\sum}\limits}
\renewcommand\Re{\mbox{Re\,}}

\renewcommand{\le}{\leqslant}

\def\la{\lambda}

\def\eps{\varepsilon}

\def\F{{\mathcal F}}

\def\bR{{\mathbb R}}

\newlength{\lenun}
\newlength{\lendu}
\begin{document}

 \begin{center}\Large{
{\bf On the Interpolation of Analytic Maps
}}\end{center}

\begin{center}\large{
A.M.Savchuk and A.A.Shkalikov}
\end{center}

\begin{center}\large{
June 22, 2013}
\end{center}

\vspace{0,3cm}

The interpolation theorem for nonlinear maps was essentially used in our paper
 \cite{SS1}, as   well, as in the papers
\cite{SS2} and \cite{SS3}. The papers dealt with nonlinear analytic maps
associated with inverse  Sturm-Liouville problems.  The required estimates
 for these maps firstly were obtained for integer  indices, and then for all intermediate values
  by nonlinear interpolation. Working on the paper \cite{SS1}  we were not aware on nonlinear
  interpolation theorems, but  assuming analyticity as additional  assumption  we succeeded
  to get the following result.

{\bf Theorem 1.} {\sl
Let $(E_0,E_1)$ and $(H_0,H_1)$   be a pair of Banach spaces, such that
 $E_1$ is densely embedded into  $E_0$ and $\|x\|_{E_0}\le
\|x\|_{E_1}$, while  $H_1$ is densely and continuously embedded into $H_0$.
Denote by
$E_\theta=[E_0,E_1]_\theta$ and $H_\theta=[H_0,H_1]_\theta$  the spaces
obtained by the method of complex interpolation for
$\theta\in[0,1]$. Denote also by $B_\theta(0,R)$ the ball  of radius $R$ centered
at zero in the space $E_\theta$. Let $\Phi$  be an analytic map acting from the ball $B_0(0,R)$
  into the space  $H_0$, such that
\begin{equation*}
\|\Phi(x)\|_{H_0}\le
C_0(R)\|x\|_{E_0}, \qquad x\in B_0(0,R).
\end{equation*}
Let $\Phi$ maps also the ball
$B_1(0,R)$ into $H_1$, and
\begin{equation}\label{M}
\|\Phi(x)\|_{H_1}\le C_1(R)\|x\|_{E_1},  \qquad  x\in B_1(0,R).
\end{equation}
 Then for all  $\theta\in [0,1]$ the map  $\Phi$ maps the ball $B_\theta(0,r)$ in the space
  $E_\theta$ of radius  $r\in(0,R)$ into
$H_\theta$, and
$$
\|\Phi(x)\|_{H_\theta}\le
C_0^{1-\theta}C_1^\theta\frac{R}{R-r}\|x\|_{E_\theta}, \qquad x\in B_\theta(0,r).
$$
}

This theorem is sufficient to complete the proof of the main Theorem 6.1 of the paper  \cite {SS1}
(see \S 6 of this paper). However, soon we found a theorem in the paper of Tartar \cite{T},
which we assumed to be more general rather then the above result  (up to constants having no essential role
 in the estimates). We used Tartar theorem  in the following form.

{\bf Theorem} (Tartar \cite{Ta}).
{\sl Let $(E_0,E_1)$ and $(H_0,H_1)$ be two  pairs of Banach spaces with dense and continuous  embeddings
 $E_1\hookrightarrow E_0$,  $H_1\hookrightarrow H_0$.
Let $\Phi$  be a nonlinear map from $E_0$ into
$E_1$, which maps  $E_1$ into $H_1$  and satisfies the following conditions:
There are positive increasing functions  $C_0(R)$ and $C_1(R)$, such that
\begin{equation}\label{1}
\|\Phi(\sigma) -\Phi(\tilde\sigma)\|_{H_0} \leqslant C_0(R) \|\sigma -\tilde\sigma\|_{E_0},
\end{equation}
if $\max\{\|\sigma\|_{E_0}, \|\tilde\sigma\|_{E_0}\} \leqslant R$, and
\begin{equation}\label{2}
\|\Phi(\sigma)\|_{H_1} \leqslant C_1(R) \|\sigma\|_{E_1}.
\end{equation}
Then $\Phi$ maps the intermediate spaces  $[E_0, E_1]_\theta $
into $H_\tau = [H_0, H_1]_\theta$  for all  $0\leqslant \theta \leqslant 1$,  and
$$
\|\Phi (\sigma)\|_{H_\theta} \leqslant C_\theta(R) \|\sigma\|_{E_\theta},
$$
where $C_\theta(R)$ is an increasing function on $R$.
}

We understood this theorem not correctly. We assumed that the assertion of Theorem  is valid provided
that estimate \eqref{1} holds for  $\max\{\|\sigma\|_{E_0}, \|\tilde\sigma\|_{E_0}\} \leqslant R$,
while estimate  \eqref{2} holds for $\|\sigma\|_{E_1} \leqslant R$. However, it is assumed in
 Tartar Theorem that estimate  \eqref{2}
also holds for  $\|\sigma\|_{E_0}\leqslant R$. The last assumption is too strong and in such a
form Tartar Theorem becomes useless for our purposes.  Prof. T.Kappeler paid our attention
to this mistake and we are  are very grateful to him for this note.
 At the same time we informed him that we
 had an independent approach to nonlinear interpolation and possessed
 the proof of the above Theorem 1 where the required estimate was obtained (in contrast to Tartar Theorem)
 under additional assumption of analyticity.

The goal of this  paper is to present the proof of Theorem 1.
We note that the definition of analytic  maps in Banach spaces and basic facts  related to this topic
 can be found in the book \cite{PT}. The basic facts on the interpolation in Banach spaces can be found in
 books   \cite{BL}  and \cite{T}, for example.

First we prove the following lemma.

{\bf Lemma }. {\sl
Let  $(E_0,E_1)$ and  $(H_0,H_1)$ be two pairs of Banach spaces with dense and continuous embeddings
$E_1\hookrightarrow E_0$ and $H_1\hookrightarrow H_0$.
Let  $E_\theta=[E_0,E_1]_\theta$ and
$H_\theta=[H_0,H_1]_\theta$  be intermediate  spaces obtained by the method of
complex interpolation for  $\theta\in[0,1]$. Assume that
$\Phi$ is a homogeneous map of degree  $k$,  i.e.
$\Phi(\la x)=\la^k \Phi(x)$,  moreover,  $\Phi$ maps  $E_0$ into
$H_0$ and  $E_1$ into $H_1$.  Assume also that
$$
\|\Phi(x)\|_{H_0}\le
M_0\|x\|_{E_0}^k
$$
and
$$
\|\Phi(x)\|_{H_1}\le M_1\|x\|_{E_1}^k.
$$
Then $\Phi$ maps the space  $E_\theta$ into $H_\theta$ and
$$
\|\Phi(x)\|_{H_\theta}\le
M_0^{1-\theta}M_1^\theta\|x\|_{E_\theta}^k.
$$
}
{\bf Proof}.
First, let us remind  some facts from the theory of complex interpolation
in Banach spaces. Let  $f$ be a bounded  analytic function in the strip
 $\Re z\in(0,1)$ with values in the space $E_0$ and continuous in the closed
 strip  $\Re z\in[0,1]$. Moreover, assume that the values  $f(it)$ for $t\in\bR$ lie in
 $E_0$, but the values  $f(1+it)$ for $t\in\bR$ lie in $E_1$. Let also the functions
 $f(it)$ and  $f(1+it)$ be continuous in the spaces $E_0$ and $E_1$, respectively, and tend
 to zero in the corresponding norms as  $|t|\to\infty$. Denote by $\F=\F(E_0,E_1)$
 the space of such functions endowed with the norm
$$
\|f\|_{\F}:=\max\{\max\limits_{t\in\bR}\|f(it)\|_{E_0},\max\limits_{t\in\bR}\|f(1+it)\|_{E_1}\}.
$$
Then the space  $\F$ is complete and  the interpolation space
$E_\theta$ is defined by the set
$E_\theta=\{f(\theta):\,f\in\F\}$ endowed with the norm
$$
\|x\|_{E_\theta}:=\inf\{\|f\|_{\F}:\,f\in\F,\ f(\theta)=x\}.
$$
Now, let us start proving lemma. For any function  $f\in\F$
define the function  $g_f=g(z):=M_0^{z-1}M_1^{-z}\Phi(f(z))$ and notice that
 $g$ is analytic in the strip  $\Re z\in(0,1)$ and continuous in the closed strip
  $\Re z\in[0,1]$ as the function with values  in $H_0$. For  $z=it$ we have
   $\|g(it)\|_{H_0}\le \|f(it)\|^k_{E_0}$, therefore  $\|g(it)\|_{H_0}\to0$ as
    $|t|\to\infty$.  For  $z=1+it$ it takes the values in the space  $H_1$, and due to  the
    inequality $\|g(1+it)\|_{H_1}\le \|f(1+it)\|^k_{E_1}$ we have
$\|g(1+it)\|_{H_1}\to0$ as $|t|\to\infty$. This implies that
$g\in\F(H_0,H_1)$ and  $\|g\|_{\F(H_0,H_1)}\le\|f\|^k_{\F(E_0,E_1)}$.
For given $x\in E_\theta$ and small $\eps>0$ take a function $f\in\F(E_0,E_1)$,
such that
$$
\|x\|_{E_\theta}\le\|f\|_{\F(E_0,E_1)}\le\|x\|_{E_\theta}+\eps.
$$
Then $\Phi(x)=\Phi(f(\theta))=M_0^{1-\theta}M_1^\theta
g_f(\theta)$, hence
$$
\|\Phi(x)\|_{H_\theta}\le
M_0^{1-\theta}M_1^\theta\|g_f\|_{\F(H_0,H_1)}\le
M_0^{1-\theta}M_1^\theta\|f\|^k_{\F(E_0,E_1)}\le
M_0^{1-\theta}M_1^\theta\left(\|x\|_{E_\theta}+\eps\right)^k.
$$
This gives the  desired result, since the number  $\eps>0$ can be chosen arbitrary. Lemma is proved.
.
\medskip
Now, let us pass to the proof of Theorem 1. Represent the map $\Phi$ in the ball
 $B_0(0,R)$ as the series
 $$
 \Phi(h)=\sum_{n=0}^\infty P_n(h)
 $$
  (see \cite{PT}, for example). Here  $P_n$ is a homogeneous map of the degree $n$,
  acting from the space  $E_0$ into
$H_0$ and defined by the Cauchy formula
$$
P_0=0, \qquad P_n(h)=\frac1{2\pi i}\int_{|\xi|=\rho}\frac{\Phi(\xi
h)}{\xi^{n+1}}d\xi,
$$
where the integral does not depend on  $\rho\in(0,R/\|h\|)$. Tending  $\rho \to R/\|h\|$ we obtain
the estimate $\|P_n(h)\|_{H_0}\le \ C_0(R)R^{1-n}\|h\|^n_{E_0}$.
The same expansion is valid for the map $\Phi$  in the ball $B_1(0,R)$.
By assumption of Theorem this ball is embedded into the ball
$B_0(0,R)$, therefore  the maps  $P_n$ in the space  $E_1$  are the restrictions onto $E_1$ of
the maps $P_n$ in the space  $E_0$.  Hence, we can apply the same arguments and obtain the estimate
$\|P_n(h)\|_{H_1}\le \ C_1(R)R^{1-n}\|h\|^n_{E_1}$.By virtue of Lemma 2   $P_n$ maps  elements  $h\in
E_\theta$ into the space  $H_\theta$, moreover,
$$
\|P_n(h)\|_{H_\theta}\le C_0(R)^{1-\theta}C_1(R)^\theta
R^{1-n}\|h\|^n_{E_\theta}.
$$
Then the map $\Phi$ defined in the ball $B_\theta(0,r)$, $r<R$, by the series
$$
\Phi(h)=\sum\limits_{n=0}^\infty P_n(h)
$$
maps this ball into the space  $H_\theta$, and
\begin{multline*}
\|\Phi(h)\|_{H_\theta}\le
C_0(R)^{1-\theta}C_1(R)^\theta\suml_{n=1}^\infty
R^{1-n}\|h\|^n_{E_\theta}\le
\\
\le
C_0(R)^{1-\theta}C_1(R)^\theta\|h\|_{E_\theta}\frac{R}{R-\|h\|_{E_\theta}}\le
C_0^{1-\theta}C_1^\theta\frac{R}{R-r}\|h\|_{E_\theta}.
\end{multline*}
This ends the proof of Theorem 1.

{\bf Remark}.  The estimate  $\|x\|_{E_0}\leqslant \|x\|_{E_1}$ in Theorem 1 can be replaced
by the estimate
 $\|x\|_{E_0}\leqslant   C\|x\|_{E_1}$, i.e. by the assumption that the embedding
  $E_1\hookrightarrow E_0$ is continuous. Assuming the last assumption one can achieve the first estimate
  by passing in one of the spaces to equivalent norm. However, the constants in \eqref{M} have to be changed
  if we return to the original norm.



The work is supported by Russian Fund of Fundamental Research, grant number 13-01-00705.

{\it Acknowledgment}.  We thank Professors T.Kappeler and
P.Topalov who initiated the publication of this  paper. They
informed us also \cite{KT} that they found an independent proof of
another version of Theorem 1 for the case when $E_\theta$ are
Sobolev spaces of periodic functions and $H_\theta$ are weighed
spaces of complex sequences.

\bigskip

A.M.Savchuk,\\
Lomonosov Moscow State University,\\
Department of Mechanics and Mathematics \\
artem\_savchuk@mail.ru

A.A.Shkalikov,\\
Lomonosov Moscow State University,\\
Department of Mechanics and Mathematics\\
 ashkalikov@yahoo.com
\end{document}